\documentstyle[11pt]{article}
\def\b#1{{\bf #1}}
\def\i#1{{\it #1}}

\begin{document}
\title{Some Formulae for Bivariate Fibonacci and Lucas Polynomials}
\author{Mario Catalani\\
Department of Economics, University of Torino\\ Via Po 53, 10124 Torino, Italy\\
mario.catalani@unito.it}
\date{}
\maketitle
\begin{abstract}
\small{We derive a collection of identities for
bivariate Fibonacci and
Lucas polynomials using essentially a matrix approach as well as properties
of such polynomials when the variables $x$ and $y$ are replaced by
polynomials. A wealth of combinatorial identities can be obtained for
selected values of the variables.}
\end{abstract}

\section{Introduction}
We define bivariate Fibonacci polynomials as
$$F_n(x,\,y)=xF_{n-1}(x,\,y)+yF_{n-2}(x,\,y),\qquad F_0(x,\,y)=0,\,
F_1(x,\,y)=1,$$ and
bivariate Lucas polynomials as
$$L_n(x,\,y)=xL_{n-1}(x,\,y)+yL_{n-2}(x,\,y),\qquad L_0(x,\,y)=2,\,
L_1(x,\,y)=x.$$
We assume $x\not= 0,\,y\not= 0,\,x^2+4y\not= 0$.

\noindent
The roots of the characteristic equation are
$$\alpha\equiv\alpha(x,\,y)={x+\sqrt{x^2+4y}\over 2},\quad\quad
\beta\equiv\beta(x,\,y)={x-\sqrt{x^2+4y}\over 2}.$$
We have $\alpha+\beta=x,\,\alpha\beta=-y$ and $\alpha-\beta=\sqrt{x^2+4y}$.
The Binet's forms are
$$F_n(x,\,y)={\alpha^n-\beta^n\over \alpha-\beta},
\quad\quad L_n(x,\,y)=\alpha^n+\beta^n.$$
Many basic facts concerning these kinds of polynomials can be found in
\cite{mario7}.

\section{Some Matrix Techniques}
Let
$$\b{B}=\left [\begin{array}{cc} x^2+2y&x\\xy&2y\end{array}\right ].$$
Then
$$\b{B}=2y\b{I}+x\b{A},$$
where
$$\b{A}=\left [\begin{array}{cc} x&1\\y&0\end{array}\right ].$$
It follows
\begin{equation}
\label{eq:one}
\b{B}^n=\sum_{k=0}^n{n\choose k}(2y)^{n-k}x^k\b{A}^k.
\end{equation}
It is known that (e.g. \cite{mario7})
$$\b{A}^n= \left [\begin{array}{cc} F_{n+1}(x,\,y)&F_n(x,\,y)\\
yF_n(x,\,y)&yF_{n-1}(x,\,y)\end{array}\right ].$$
Now using the method exposed in \cite{mclaughlin} to $\b{B}$ we
find, after some manipulations, that the element $(1,\,2)$ of $\b{B}^n$
is given by
$$x\sum_{k=0}^{\left\lfloor{n-1\over 2}\right\rfloor}
{n-1-k\choose k}(x^2+4y)^{n-1-k}(-y)^k.$$
The corresponding element of the RHS of Equation~\ref{eq:one} is
$$\sum_{k=0}^n{n\choose k}(2y)^{n-k}x^kF_k(x,\,y).$$
Then, using the Binet's form,
$$x\sum_{k=0}^{\left\lfloor{n-1\over 2}\right\rfloor}
{n-1-k\choose k}(x^2+4y)^{n-1-k}(-y)^k
={1\over \alpha-\beta}\sum_{k=0}^n{n\choose k}(2y)^{n-k}x^k
(\alpha^k-\beta^k).$$
Now
\begin{eqnarray*}
\sum_{k=0}^n{n\choose k}(2y)^{n-k}x^k
\alpha^k&=&(2y+x\alpha)^n\\
&=&\alpha^n(\alpha-\beta)^n.
\end{eqnarray*}
Analogously for the $\beta$ part. After replacing $n$ by $2n$ we get the
identity
\begin{equation}
\label{eq:prima}
x\sum_{k=0}^{n-1}{2n-1-k\choose k}(x^2+4y)^{n-k-1}(-y)^k=
F_{2n}(x,\,y).
\end{equation}
Now
$$\b{BA}=\left [\begin{array}{cc}x^3+3xy&x^2+2y\\x^2y+2y^2&
xy\end{array}\right ],$$
and
$${\rm tr}(\b{BA})=x^3+4xy=x(x^2+4y),$$
$$\vert \b{BA}\vert=\vert\b{B}\vert\,\vert\b{A}\vert=
-y(x^2y+4y^2)=-y^2(x^2+4y).$$
Applying the previous method this time to $\b{BA}$  we get that
the element $(1,\,2)$ of $(\b{BA})^n$ is given by
\begin{eqnarray*}
&&(x^2+2y)\sum_{k=0}^{\left\lfloor n-1\over 2\right\rfloor}
{n-1-k\choose k}x^{n-1-2k}(x^2+4y)^{n-1-2k}y^{2k}(x^2+4y)^k=\\
&&\quad\quad =
(x^2+2y)\sum_{k=0}^{\left\lfloor n-1\over 2\right\rfloor}
{n-1-k\choose k}x^{n-1-2k}(x^2+4y)^{n-1-k}y^{2k}.
\end{eqnarray*}
On the other hand (see \cite{mario7})
$$\b{BA}=(2y\b{I}+x\b{A})\b{A}=2y\b{A}+x\b{A}^2,$$
so that
\begin{eqnarray*}
(\b{BA})^n&=&\sum_{k=0}^n{n\choose k}(2y)^{n-k}\b{A}^{n-k}x^k\b{A}^{2k}\\
&=&\sum_{k=0}^n{n\choose k}(2y)^{n-k}x^k\b{A}^{n+k}.
\end{eqnarray*}
Since the element $(1,\,2)$ of $\b{A}^n$ is $F_n(x,\,y)$, we have that
the element $(1,\,2)$ of $(\b{BA})^n$ is
$$\sum_{k=0}^n{n\choose k}(2y)^{n-k}x^kF_{n+k}(x,\,y)=
{1\over\alpha-\beta}\sum_{k=0}^n{n\choose k}(2y)^{n-k}x^k
(\alpha^{n+k}-\beta^{n+k}).$$
Now
\begin{eqnarray*}
{1\over\alpha-\beta}\sum_{k=0}^n{n\choose k}(2y)^{n-k}x^k
\alpha^{n+k}&=&
{\alpha^n\over\alpha-\beta}\sum_{k=0}^n{n\choose k}(2y)^{n-k}(x
\alpha)^k\\
&=&
{\alpha^n\over\alpha-\beta}(2y+
x\alpha)^n\\
&=&{\alpha^n\over\alpha-\beta}\alpha^n(\alpha-\beta)^n\\
&=&{\alpha^{2n}(\alpha-\beta)^n\over\alpha-\beta}.
\end{eqnarray*}
In the same way
\begin{eqnarray*}
{1\over\alpha-\beta}\sum_{k=0}^n{n\choose k}(2y)^{n-k}x^k
\alpha^{n+k}&=&
{\beta^n\over\alpha-\beta}\sum_{k=0}^n{n\choose k}(2y)^{n-k}(x
\beta)^k\\
&=&
{\beta^n\over\alpha-\beta}(2y+
x\beta)^n\\
&=&{\beta^n\over\alpha-\beta}\beta^n(-1)^n(\alpha-\beta)^n\\
&=&{\beta^{2n}(-1)^n(\alpha-\beta)^n\over\alpha-\beta}.
\end{eqnarray*}
Subtracting we get
$$(\alpha-\beta)^n\left ({\alpha^{2n}\over\alpha-\beta}-
{\beta^{2n}(-1)^n\over\alpha-\beta}\right ).$$

\noindent
Replacing $n$ with $2n$ and equating the two expressions for the element
$(1,\,2)$ we get
\begin{eqnarray*}
&&(x^2+2y)\sum_{k=0}^{n-1}
{2n-1-k\choose k}x^{2n-1-2k}(x^2+4y)^{2n-1-k}y^{2k}\\
&&\qquad\qquad =(x^2+4y)^nF_{4n}(x,\,y),
\end{eqnarray*}
that is
\begin{equation}
\label{eq:unouno}
(x^2+2y)\sum_{k=0}^{n-1}
{2n-1-k\choose k}x^{2n-1-2k}(x^2+4y)^{n-1-k}y^{2k}
=F_{4n}(x,\,y).
\end{equation}
Using Equation~\ref{eq:prima}
and replacing here $n$ by $2n$ we obtain the identity
\begin{eqnarray*}
&&
(x^2+2y)\sum_{k=0}^{n-1}
{2n-1-k\choose k}x^{2n-1-2k}(x^2+4y)^{n-1-k}y^{2k}\\
&&\quad\quad =
x\sum_{k=0}^{2n-1}
{4n-1-k\choose k}(x^2+4y)^{2n-1-k}(-y)^k.
\end{eqnarray*}

\noindent
\section{Polynomials as Arguments}
First of all, from
$$F_n(x,\,y)={\alpha^n-\beta^n\over\alpha-\beta},\qquad
L_n(x,\,y)=\alpha^n+\beta^n,$$
solving for $\alpha$ and $\beta$ we get
\begin{equation}
\label{eq:uno}
\alpha^n={L_n(x,\,y)+(\alpha-\beta)F_n(x,\,y)\over 2},
\end{equation}
\begin{equation}
\label{eq:due}
\beta^n={L_n(x,\,y)-(\alpha-\beta)F_n(x,\,y)\over 2}.
\end{equation}

\noindent
Secondly
\begin{eqnarray*}
L_{n}^2(x,\,y)&=&\alpha^{2n}+\beta^{2n}+2(-1)^ny^n\\
&=&\alpha^{2n}+\beta^{2n}-2(-1)^{n+1}y^n,
\end{eqnarray*}
\begin{eqnarray*}
(\alpha-\beta)^2F_{n}^2(x,\,y)&=&
\alpha^{2n}+\beta^{2n}-2(-1)^ny^n\\
&=&\alpha^{2n}+\beta^{2n}+2(-1)^{n+1}y^2.
\end{eqnarray*}
It follows
\begin{equation}
\label{eq:sei}
L_{n}^2(x,\,y)+(-1)^{n+1}4y^{n}=(x^2+4y)F_{n}^2(x,\,y).
\end{equation}
Then
\begin{eqnarray*}
\alpha\left (L_{n}(x,\,y),\,(-1)^{n+1}y^n\right )&=&
{L_{n}(x,\,y)+\sqrt{L_{n}^2(x,\,y)+4(-1)^{n+1}y^{n}}\over 2}\\
&=&{L_{n}(x,\,y)+\sqrt{x^2+4y}F_{n}(x,\,y)\over 2}\\
&=&\alpha^{n}(x,\,y).
\end{eqnarray*}
$$\beta\left (L_{n}(x,\,y),\,(-1)^{n+1}y^{n}\right )
=\beta^{n}(x,\,y).$$
Then
\begin{eqnarray}
\label{eq:cinque}
&&F_n\left (L_{k}(x,\,y),\,(-1)^{k+1}y^k\right )\nonumber\\
&&\quad\quad ={\alpha^n\left (L_{k}(x,\,y),\,(-1)^{k+1}y^k\right )-
\beta^n\left (L_{k}(x,\,y),\,(-10^{k+1}y^k\right )\over
\alpha\left (L_{k}(x,\,y),\,(-1)^{k+1}y^k\right )-
\beta\left (L_{k}(x,\,y),\,(-1)^{k+1}y^k\right ) }\nonumber\\
&&\quad\quad =
{\alpha^{nk}(x,\,y)-
\beta^{nk}(x,\,y)\over
\alpha^{k}(x,\,y)-
\beta^{k}(x,\,y)}\nonumber\\
&&\quad\quad ={F_{nk}(x,\,y)\over F_{k}(x,\,y)}.
\end{eqnarray}
As a consequence,
we can see that $F_k(x,\,y)$ divides $F_{kn}(x,\,y)$, $\forall \,k,\,n$.

\noindent
The analogous of Equation~\ref{eq:cinque} for Lucas polynomials is
\begin{equation}
\label{eq:lucas}
L_n\left (L_{k}(x,\,y),\,(-1)^{k+1}y^{k}\right )=L_{nk}(x,\,y).
\end{equation}
\b{Some examples.}
$$F_{2n}(x,\,y)=F_{2}(x,\,y)F_{n}(L_{2}(x,\,y),\,-y^2)
=xF_{n}(x^2+2y,\,-y^2),$$
$$
F_{3n}(x,\,y)=F_{3}(x,\,y)F_{n}(L_{3}(x,\,y),\,y^3)
=(x^2+y)F_{n}(x^3+3xy,\,y^3),$$
$$F_{4n}(x,\,y)=F_{4}(x,\,y)F_{n}(L_{4}(x,\,y),\,-y^4)
=x(x^2+2y)F_{n}(x^4+4x^2y+2y^2,\,-y^4).$$
\bigskip

\noindent
Equation~\ref{eq:sei} and Equation~\ref{eq:lucas}
can be used to derive many other identities. For
instance the \i{Simpson} formula (see \cite{mario7})
\begin{equation}
\label{eq:simpson0}
L_n(x,\,y)L_{n+2}(x,\,y)-L_{n+1}^2(x,\,y)=
(-1)^ny^n(x^2+4y)
\end{equation}
setting $x=L_k(x,\,y),\,y=(-1)^{k+1}y^k$ becomes
\begin{eqnarray*}
&&L_n\left (L_k(x,\,y),\,(-1)^{k+1}y^k\right )
L_{n+2}\left (L_k(x,\,y),\,(-1)^{k+1}y^k\right )+\\
&&\qquad\qquad\quad -L_{n+1}^2\left (L_k(x,\,y),\,(-1)^{k+1}y^k\right )
\\&&\quad
=(-1)^n((-1)^{k+1}y^k)^n(L_k^2(x,\,y)+4(-1)^{k+1}y^k)\\
&&\quad = (-y)^{nk}(x^2+4y)F_k^2(x,\,y),
\end{eqnarray*}
that is
$$L_{kn}(x,\,y)L_{k(n+2)}(x,\,y)-L_{k(n+1)}^2(x,\,y)
=(-y)^{nk}(x^2+4y)F_k^2(x,\,y).$$
Equation~\ref{eq:prima}
with $x=L_k(x,\,y),\,y=(-1)^{k+1}y^k$ becomes
\begin{eqnarray*}
&&L_k(x,\,y)\sum_{r=0}^{n-1}{2n-1-r\choose r}(x^2+4y)^{n-1-r}
F_k^{2(n-1-r)}(x,\,y)(-y)^{rk}\\
&&\quad\qquad =F_{2n}(L_k(x,\,y),\,(-1)^{k+1}y^k).
\end{eqnarray*}
Now using the identity
$$L_n^2(x,\,y)+2(-1)^{n+1}y^n=L_{2n}(x,\,y),$$
we have
\begin{eqnarray*}
&&F_{2n}\left (L_k(x,\,y),\,(-1)^{k+1}y^k\right )\\
&&\quad\quad =L_k(x,\,y)
F_n\left (L_k^2(x,\,y)+ 2(-1)^{k+1}y^k,\,-((-1)^{k+1}y^k)^2\right )\\
&&\quad\quad
=L_k(x,\,y)
F_n\left (L_{2k}(x,\,y),\,-y^{2k}\right ).
\end{eqnarray*}
Finally
$$
\sum_{r=0}^{n-1}{2n-1-r\choose r}(x^2+4y)^{n-1-r}
F_k^{2(n-1-r)}(x,\,y)(-y)^{rk}=
F_n\left (L_{2k}(x,\,y),\,-y^{2k}\right ).$$
Then, since $-y^{2k}=(-1)^{2k+1}y^{2k}$, Equation~\ref{eq:cinque} allows to
write
$$
\sum_{r=0}^{n-1}{2n-1-r\choose r}(x^2+4y)^{n-1-r}
F_k^{2(n-1-r)}(x,\,y)(-y)^{rk}=
{F_{2kn}(x,\,y)\over F_{2k}(x,\,y)}.$$
Equation~\ref{eq:cinque} and Equation~\ref{eq:lucas} allow to derive new
identities. For example, the identity
\begin{equation}
\label{eq:simpson1}
yF_{n-1}(x,\,y)+F_{n+1}(x,\,y)=L_{n}(x,\,y)
\end{equation}
upon substitution of $x$ with $L_{k}(x,\,y)$ and $y$ with $(-1)^{k+1}y^k$
becomes
$$(-1)^{k+1}y^k{F_{k(n-1)}(x,\,y)\over F_{k}(x,\,y)}+
{F_{k(n+1)}(x,\,y)\over
F_{k}(x,\,y)}=L_{nk}(x,\,y),$$
that is
$$(-1)^{k+1}y^kF_{k(n-1)}(x,\,y)+
F_{k(n+1)}(x,\,y)
=F_{k}(x,\,y)L_{nk}(x,\,y).$$
Or the identity
$$L_{n+2}^2(x,\,y)+yL_{n+1}^2(x,\,y)=(x^2+2y)
L_{2n+2}(x,\,y)+xyL_{2n+1}(x,\,y)$$
that becomes
\begin{eqnarray*}
&&L_{k(n+2)}^2(x,\,y)+(-1)^{k+1}y^kL_{k(n+1)}^2(x,\,y)\\
&&\quad\quad =L_{2k}(x,\,y)
L_{k(2n+2)}(x,\,y)+(-1)^{k+1}y^kL_{k}(x,\,y)L_{k(2n+1)}(x,\,y).
\end{eqnarray*}
Equation~\ref{eq:sei} can be rewritten as
$$L_n^2(x,\,y)=(x^2+4y)F_n^2(x,\,y)+(-1)^n4y^n.$$
Then
$$\alpha\left (\sqrt{x^2+4y}F_n(x,\,y),\,(-1)^ny^n\right )=\alpha^n(x,\,y),$$
$$\beta\left (\sqrt{x^2+4y}F_n(x,\,y),\,(-1)^ny^n\right )=-\beta^n(x,\,y).$$
It follows
\begin{equation}
\label{eq:simpson2}
F_{2n+1}\left (\sqrt{x^2+4y}F_k(x,\,y),\,(-1)^ky^k\right )=
{L_{k(2n+1)}(x,\,y)\over L_{k}(x,\,y)},
\end{equation}
\begin{equation}
\label{eq:simpson3}
F_{2n}\left (\sqrt{x^2+4y}F_k(x,\,y),\,(-1)^ky^k\right )=
(\alpha-\beta){F_{2kn}(x,\,y)\over L_{k}(x,\,y)}.
\end{equation}
For the Lucas polynomials we have
\begin{equation}
\label{eq:simpson4}
L_{2n+1}\left (\sqrt{x^2+4y}F_k(x,\,y),\,(-1)^ky^k\right )=
(\alpha-\beta)F_{k(2n+1)}(x,\,y),
\end{equation}
\begin{equation}
\label{eq:simpson5}
L_{2n}\left (\sqrt{x^2+4y}F_k(x,\,y),\,(-1)^ky^k\right )=
L_{2kn}(x,\,y).
\end{equation}
Now identity~\ref{eq:simpson1} replacing $x$ with $\sqrt{x^2+4y}F_k(x,\,y)$
and $y$ with $(-1)^ky^k$ becomes, using Equations~\ref{eq:simpson2},
\ref{eq:simpson3},~\ref{eq:simpson4},
\ref{eq:simpson5},
$$(-1)^ky^kL_{k(2n-1)}(x,\,y)+L_{k(2n+1)}(x,\,y)=
L_{2kn}(x,\,y) L_{k}(x,\,y),$$
and
$$(-1)^ky^kF_{2kn}(x,\,y)+F_{k(2n+2)}(x,\,y)=
F_{k(2n+1)}(x,\,y) L_{k}(x,\,y).$$
In an analogous way the Simpson formula (Equation~\ref{eq:simpson0})
becomes
$$L_{2kn}(x,\,y)L_{k(2n+2)}(x,\,y)-(\alpha-\beta)^2F_{k(2n+1)}^2(x,\,y)
=y^{2nk}L_{k}^2(x,\,y),$$
and
$$(\alpha-\beta)^2F_{k(2n-1)}(x,\,y)F_{k(2n+1)}(x,\,y)-
L_{2kn}^2(x,\,y)
=-y^{k(2n-1)}L_{k}^2(x,\,y).$$

\end{document}